\documentclass[11pt,a4paper]{amsart}
\usepackage{amssymb,amsmath}
\usepackage[hidelinks]{hyperref}
\newtheorem{theorem}{Theorem}
\newtheorem{lemma}[theorem]{Lemma}
\newtheorem{corollary}[theorem]{Corollary}

\begin{document}

\title[Automatic continuity of biseparating maps]{Automatic continuity of linear maps preserving zero products or range orthogonality in both directions}

\author{Ngai-Ching Wong}

\address[Ngai-Ching Wong]{Department of Applied Mathematics, National Sun Yat-sen University,
Kaohsiung 80424, Taiwan.}
\email{wong@math.nsysu.edu.tw}

\begin{abstract}
\emergencystretch=1em
Let $T$ be a complex-linear surjection from a $C^*$-algebra onto a Banach algebra $B$.
If $T$ preserves zero products in both directions, then $B$ is semisimple and $T$ is automatically bounded.
When $B$ is a $C^*$-algebra, $T$ is an algebra isomorphism multiplied by an invertible central multiplier.
For a complex-linear surjection between $C^*$-algebras preserving range orthogonality in both directions, we also obtain automatic boundedness; in this case the map is a $*$-isomorphism followed by multiplication on the right by an invertible right multiplier.
Neither algebra need be unital.
\end{abstract}


\thanks{This research is partially supported by the Taiwan NSTC grant 114-2115-M-110-002-MY2.}

\date{November 26, 2026}

\maketitle

\section{Introduction}

For algebras $A,B$, call a complex-linear map $T:A\to B$
\emph{biseparating} if
\begin{align}\label{eq:bisep}
 ab=0\quad\Longleftrightarrow\quad T(a)T(b)=0\qquad(a,b\in A).
\end{align}
Araujo and Jarosz asked whether bijective biseparating maps between arbitrary
$C^*$-algebras are continuous and are central multiples of algebra
isomorphisms \cite[Conjecture~1]{AJ}.  Positive answers are known for abelian $C^*$-algebras \cite{Jarosz90, JW96}, CCR
algebras with Hausdorff spectrum \cite{LW}, von Neumann algebras \cite{LTW}, and $AW^*$-algebras
\cite{LCLW}.

A related question asks whether a complex-linear surjection between $C^*$-algebras is automatically bounded if it preserves range orthogonality in both directions:
\begin{equation}\label{eq:range}
a^*b=0\quad\Longleftrightarrow\quad T(a)^*T(b)=0\qquad(a,b\in A).
\end{equation}
Positive answers were obtained for continuous-trace $C^*$-algebras \cite{Tsai}, von Neumann algebras \cite{LTW}, and $AW^*$-algebras \cite{LCLW}. For algebraic results on preservers of range orthogonality, see \cite{Hamhalter}.

We answer the conjecture of Araujo and Jarosz and settle the range orthogonality case of \cite[Conjecture~3.9]{survey} in Theorem \ref{thm:transport}.
Both proofs use Brown's theorem on local maps to control multiplication operators transported by $T$.
For zero product preservers the range may even be an arbitrary Banach algebra: its semisimplicity follows from the preservation assumption. For $C^*$-algebra ranges we obtain the two factorizations in Corollary~\ref{cor:factor}.

For surjective maps, preservation in one direction does not suffice. Jarosz \cite{Jarosz90} constructed a discontinuous 
 linear surjection from the algebra $c$ of convergent sequences onto $\ell^\infty$ preserving 
 pairs with zero products. Since both algebras are commutative, this map also preserves range orthogonality in one direction.

\section{Results}

\begin{lemma}\label{lem:local}
Let $A$ be a $C^*$-algebra and let $F:A\to A$ be a complex linear map.
\begin{enumerate}
  \item[(a)] If $xy=0$ implies $F(x)y=0$ for all $x,y\in A$, then $F$ is bounded and $F(xa)=F(x)a$.
  \item[(b)] If $yx=0$ implies $yF(x)=0$ for all $x,y\in A$, then $F$ is bounded and $F(ax)=aF(x)$.
\end{enumerate}
\end{lemma}
\begin{proof}
We check locality in the precise sense of \cite[Definition~1.6]{Brown}.
Let $q$ be an open projection of $A^{**}$ and suppose that $xq=0$.
Take a positive contractive approximate identity $(e_\lambda)$ for
the hereditary subalgebra $qA^{**}q\cap A$. Since $xe_\lambda=0$,
we have $F(x)e_\lambda=0$ for every $\lambda$. The net
$(e_\lambda)$ converges weak$^*$ to $q$ in $A^{**}$, so
$F(x)q=0$. Thus $F$ is local, and Brown's
Corollary~5.4 \cite{Brown}, with both open projections equal to $1$,
gives boundedness, and Proposition~5.5 there gives the module identity. For the other assertion, apply the first to
$x\mapsto F(x^*)^*$.
\end{proof}

\begin{lemma}\label{lem:transport}
Let $U:X\to D$ be a linear bijection from a Banach space onto a semisimple
Banach algebra. Write $L_d(y)=dy$ and $R_d(y)=yd$. If
$U^{-1}L_dU$ and $U^{-1}R_dU$ are bounded on $X$ for every $d\in D$,
then $U$ and $U^{-1}$ are bounded.
\end{lemma}
\begin{proof}
Give $X$ the product $x\star y=U^{-1}(U(x)U(y))$.
For $x,y\in X$, let
\[
   \lambda_x(z)=x\star z,\qquad \rho_y(z)=z\star y
   \qquad(z\in X).
\]
By hypothesis, $\lambda_x=U^{-1}L_{U(x)}U$ and
$\rho_y=U^{-1}R_{U(y)}U$ are bounded operators on $X$.
It follows
that the bilinear product $\star$ is separately continuous.

Consider the family $\{\lambda_x:\|x\|\leq 1\}$. For each fixed
$y\in X$, we have
\[
   \sup_{\|x\|\leq 1}\|\lambda_x(y)\|
   =\sup_{\|x\|\leq 1}\|\rho_y(x)\|
   =\|\rho_y\|<\infty.
\]
Thus this family is pointwise bounded on the Banach space $X$.
By the uniform boundedness principle,
\[
   c:=\sup_{\|x\|\leq 1}\|\lambda_x\|<\infty.
\]
For $x\neq 0$, bilinearity gives
$\lambda_x=\|x\|\lambda_{x/\|x\|}$, so
\[
   \|x\star y\|=\|\lambda_x(y)\|
   \leq c\|x\|\|y\|.
\]
The inequality also holds for $x=0$.
If $X=0$, the assertion is immediate. Otherwise $c>0$, since a nonzero semisimple algebra cannot have identically zero multiplication.
Thus $(X,\star)$ with norm $c\|\cdot\|$ is a Banach algebra, and $U$ is an
algebra isomorphism onto $D$. Johnson's uniqueness of complete algebra norm
topology \cite{Johnson} gives the boundedness of  $U$. Its inverse is bounded by the open mapping theorem.
\end{proof}

\begin{theorem}\label{thm:transport}
Let $T:A\to B$ be a complex-linear surjection from a $C^*$-algebra
onto a Banach algebra $B$. In either of the following cases,
$T$ is bijective and both $T$ and $T^{-1}$ are bounded:
\begin{enumerate}
\item[(a)] $T$ satisfies \eqref{eq:bisep} (in this case $B$ is necessarily semisimple);
\item[(b)] $B$ is a $C^*$-algebra and
$T$ satisfies \eqref{eq:range}.
\end{enumerate}
\end{theorem}
\begin{proof}
If $T(a)=0$, then in case (a) we have $T(a^*)T(a)=0$, and in case (b) we have $T(a)^*T(a)=0$. The reverse implication in \eqref{eq:bisep} or \eqref{eq:range}, respectively, gives $a^*a=0$. Thus $T$ is injective.

Fix $d\in B$ and put $P_d=T^{-1}L_dT$, $Q_d=T^{-1}R_dT$.
We first study these two operators without assuming that $B$ is semisimple.

In case (a), suppose that $xy=0$. By the forward
implication in \eqref{eq:bisep}, we have $T(x)T(y)=0$. Since
$P_d(x)=T^{-1}(dT(x))$, associativity in $B$ gives
\[
   T(P_d(x))T(y)=(dT(x))T(y)=d(T(x)T(y))=0.
\]
It follows from the reverse implication in \eqref{eq:bisep} that $P_d(x)y=0$. Similarly, since
$Q_d(y)=T^{-1}(T(y)d)$, we have
\[
   T(x)T(Q_d(y))=T(x)(T(y)d)=(T(x)T(y))d=0.
\]
The reverse implication in \eqref{eq:bisep} gives $xQ_d(y)=0$.
Thus
$$
xy=0 \quad\text{implies}\quad P_d(x)y=0,\ xQ_d(y)=0.
$$
The first assertion of Lemma~\ref{lem:local} gives boundedness of
$P_d$, and the second gives boundedness of $Q_d$.

It remains in case (a) to prove that $B$ is semisimple.
Define $x\star y=T^{-1}(T(x)T(y))$ on $A$. Its left and right multiplications are $P_{T(x)}$ and $Q_{T(y)}$, respectively. They are bounded by the preceding paragraph, so the uniform boundedness argument in Lemma~\ref{lem:transport} shows that $\star$ is jointly continuous in the original norm of $A$. The module identities in Lemma~\ref{lem:local} give
\[
 (ax)\star y=a(x\star y),\qquad x\star(yb)=(x\star y)b .
\]
Since $xy=0$ implies $x\star y=0$, the zero-product theorem of Alaminos, Bre\v sar, Extremera and Villena 
(\cite[Lemma 2.3]{ABEV09}, \cite[Theorem~1.1]{ABEV}) gives a bounded linear map $V:A\to A$ such that $x\star y=V(xy)$.
Every $z\in A$ is a product $xy$ by Cohen factorization. Hence the two module identities yield
\[
 V(az)=aV(z),\qquad V(zb)=V(z)b.
\]
Also, if $V(z)=0$, then $x\star y=0$, so the reverse implication in \eqref{eq:bisep} gives $z=xy=0$.

Thus $I=V(A)$ is a two-sided ideal of $A$, and $V$ is an algebra isomorphism from $(A,\star)$ onto $I$: indeed,
\[
 V(x\star y)=V(V(xy))=V(V(x)y)=V(x)V(y).
\]
Although $I$ need not be closed, its Jacobson radical is zero.
Otherwise take $0\ne z\in\operatorname{rad} I$ and put $w=z^*zz^*\in I$.
Then $wz=(z^*z)^2$ belongs to $\operatorname{rad} I$. For every nonzero scalar $\lambda$, the element $\lambda 1-wz$ is invertible in the unitization of $I$, hence in the unitization of $A$. This contradicts the fact that the nonzero positive element $(z^*z)^2$ has its norm in its spectrum. Since $B$ is algebraically isomorphic to $I$, it follows that $B$ is semisimple. Lemma~\ref{lem:transport} now applies to $T$.

In case (b), if $x^*y=0$, then $T(x)^*T(Q_d(y))=T(x)^*T(y)d=0$, whence
$x^*Q_d(y)=0$. Since $x^*$ ranges over $A$, Lemma~\ref{lem:local}(b) shows that $Q_d$ is bounded.
To treat $P_d$, we use the range orthogonality condition in the opposite direction.

Fix $a\in A$ and put $S_a=TR_aT^{-1}:B\to B$, so that
$S_a(T(x))=T(xa)$ for every $x\in A$. We first show that
$S_a(uv)=uS_a(v)$ for all $u,v\in B$.
Indeed, suppose that $u^*v=0$, and write $u=T(x)$ and $v=T(y)$.
By the reverse implication in \eqref{eq:range}, we have $x^*y=0$.
Hence $x^*(ya)=(x^*y)a=0$, and the forward implication gives
\[
   u^*S_a(v)=T(x)^*T(ya)=0.
\]
Since $u^*$ ranges over all of $B$, this says that $wv=0$ implies
$wS_a(v)=0$. The second assertion of Lemma~\ref{lem:local} now gives
$S_a(wv)=wS_a(v)$ for all $w,v\in B$. In particular,
\begin{align}\label{eq:Sa}
S_a(dv)=dS_a(v)\qquad (v\in B).
\end{align}

Recall that $S_a(T(z))=T(za)$ and $T(P_d(z))=dT(z)$ for every
$z\in A$. Set $v=T(x)$ in \eqref{eq:Sa}.
On the left,
\[
   S_a(dT(x))=S_a(T(P_d(x)))=T(P_d(x)a).
\]
On the right, we obtain
\[
   dS_a(T(x))=dT(xa)=T(P_d(xa)).
\]
The two expressions are equal. Since $T$ is injective, it follows
that $P_d(x)a=P_d(xa)$ for all $x,a\in A$.
In particular, $P_d(x)y=P_d(xy)=0$ whenever
$xy=0$. One more application of Lemma~\ref{lem:local} gives the
boundedness of $P_d$. Since the $C^*$-algebra $B$ is semisimple, Lemma~\ref{lem:transport} completes the proof.
\end{proof}

\begin{corollary}\label{cor:factor}
Let $T:A\to B$ be a complex-linear surjection between $C^*$-algebras.
\begin{enumerate}
\item[(a)] $T$ is biseparating if and only if $T(a)=h\pi(a)$,
where $h$ is an invertible central multiplier of $B$ and $\pi:A\to B$
is a bounded algebra isomorphism.
\item[(b)] $T$ preserves range orthogonality in both directions if and only if $T(a)=\pi(a)h$,
where $\pi:A\to B$ is a $*$-isomorphism and $h$ is an invertible
right multiplier of $B$.
\end{enumerate}
In each case the factors are unique, with $h=T^{**}(1_{A^{**}})$.
\end{corollary}
\begin{proof}
Theorem~\ref{thm:transport} gives boundedness and bijectivity. In case (a),
\cite[Theorem~2.5]{Wong} gives the central weight and a surjective
homomorphism, necessarily injective. In case (b), \cite[Theorem~3.6(c)]{survey}
gives invertibility of $h$ and a $*$-isomorphism $\pi(a)=T(a)h^{-1}$
onto $B$. In particular, $Bh=B$.

Conversely, in case (a) we have
$T(a)T(b)=h^2\pi(ab)$; in case (b) we have $T(a)^*T(b)=h^*\pi(a^*b)h$.
These identities give the respective equivalences. In both factorizations
$\pi^{**}$ is unital, so evaluating $T^{**}$ at the identity recovers $h$.
\end{proof}

For domain orthogonality, $ab^*=0$, apply part~(b) to the complex-linear
map $a\mapsto T(a^*)^*$. The corresponding form is $T(a)=h\pi(a)$
with an invertible left multiplier. The weight for range orthogonality
need not be central: on $M_2(\mathbb C)$, take $T(a)=a\,\mathrm{diag}(1,2)$.
Preservation of the double orthogonality $a^*b=ab^*=0$ is a different
hypothesis and is not treated here; interested readers are referred to \cite{BGP}.

\section*{Disclosure of AI Use}

The author used ChatGPT (OpenAI) to assist with exposition,
language editing, and the development of portions of the
mathematical arguments during revision. He remains
responsible for all results, proofs, references, and conclusions.

\end{document}